\documentclass[authoryear,preprint,review,12pt]{elsarticle}

\usepackage{amsmath,amssymb,amsthm,amsfonts,mathtools,dsfont,hyperref,
	xcolor,fullpage}

\hypersetup{colorlinks}

\newcommand{\B}{\mathcal B}
\newcommand{\C}{\mathcal C}
\newcommand{\E}{\mathbb E}
\newcommand{\I}{\mathcal I}
\renewcommand{\i}{\boldsymbol i}
\newcommand{\e}{\mathrm e}
\newcommand{\x}{\boldsymbol x}
\renewcommand{\d}{\mathrm d}
\renewcommand{\P}[1]{\mathbb{P}\left\{#1\right\}}
\newcommand{\PP}{\mathbb P\hs}
\renewcommand{\SS}{\mathcal S}
\newcommand{\X}{\mathbb X}
\newcommand{\Z}{\mathbb Z}

\newcommand{\hs}{\hspace{1pt}}

\makeatletter
\newcommand{\svdots}{
	\vbox{\fontsize{\sf@size}{\sf@size pt}\linespread{0.3}\selectfont
		\kern0.2\baselineskip
		\hbox{.}\hbox{.}\hbox{.}%
		\kern0.1\baselineskip
}}
\def\ps@pprintTitle{%
	\let\@oddhead\@empty
	\let\@evenhead\@empty
	\def\@oddfoot{}%
	\let\@evenfoot\@oddfoot}
\makeatother

\newtheorem{thm}{Theorem}
\newtheorem{lem}[thm]{Lemma}
\newtheorem{prop}[thm]{Proposition}

\newdefinition{rmk}{Remark}
\newproof{pf}{Proof}

\newcommand{\proofofref}{}
\newproof{proofof}{Proof of \proofofref}
\newenvironment{pot}[1]
{\renewcommand{\proofofref}{#1}\proofof}
{\endproofof}

\begin{document}
	
	\begin{frontmatter}
		
		\title{Short cycles of random permutations with cycle weights:\\point processes approach}
		
		\author[1]{Oleksii Galganov}
		\ead{galganov.oleksii@lll.kpi.ua}
		\author[1,2]{Andrii Ilienko\corref{cor1}}
		\ead{andrii.ilienko@unibe.ch}
		\affiliation[1]{organization={Igor Sikorsky Kyiv Polytechnic Institute},
			addressline={Prospect Beresteiskyi 37},
			postcode={03056},
			city={Kyiv},
			country={Ukraine}}
		\affiliation[2]{organization={Institute of Mathematical Statistics and Actuarial Science, University of Bern},
			addressline={Alpeneggstrasse 22},
			city={Bern},
			postcode={CH-3012},
			country={Switzerland}}
		\cortext[cor1]{Corresponding author}
		
		\begin{abstract}
			We study the asymptotic behavior of short cycles of 
			random permutations with cycle weights. More specifically, 
			on a specially constructed metric space whose elements encode 
			all possible cycles, we consider a point process containing all 
			information on cycles of a given random permutation on 
			$\{1,\ldots,n\}$. The main result of the paper is the 
			distributional convergence with respect to the vague topology of 
			the above processes towards a Poisson point process as $n\to\infty$ 
			for a wide range of cycle weights. As an application, we give 
			several limit theorems for various statistics of cycles.
		\end{abstract}
		
		
	\end{frontmatter}
	
	\section{Introduction}
	
	Random permutations are a classical object of combinatorial 
	probability. Uniform permutations (that is, uniformly distributed 
	random elements of the symmetric group $\SS_n$) have been studied 
	since de Montmort's matching problem. In recent years, there 
	has been an extensive literature on non-uniform permutations 
	$\sigma_n$ with cycle weights, which are defined by
	the probability distribution
	\begin{equation}\label{eq:dist}
		\P{\sigma_n=\pi} = \frac{1}{h_n n!} \prod_{k=1}^\infty \theta_k^{C_k(\pi)}, \qquad \pi \in \SS_n.
	\end{equation}
	Here $\theta_k$ are non-negative parameters, $C_k(\pi)$ stands for 
	the number of $k$-cycles in $\pi$, and $h_n$ is a normalization ensuring 
	that $\sum_{\pi\in\SS_n}\P{\sigma_n=\pi}=1$. These permutations 
	were introduced in \cite{BU11}, motivated by the theory of 
	Bose-Einstein condensate. For other applications and connections, 
	see \cite{EU14} and references therein. Note that special cases of 
	\eqref{eq:dist} are the uniform random permutation (with 
	$\theta_k=1$ for all $k$ and $h_n=1$) and Ewens random permutations 
	(with $\theta_k=\theta$ for all $k$, $h_n=\theta^{(n)}/n!$, 
	and $\theta^{(n)}$ standing for the rising factorial). 
	The latter are based on the Ewens sampling formula which 
	was introduced in population genetics and subsequently 
	found numerous applications, see \cite{C16}. A rich theory of 
	Ewens permutations was developed in \cite{ABT03}.
	
	An important subject of study of random permutations is their 
	cycle structure and, in particular, the asymptotic statistics of 
	short cycles (that is, cycles of bounded length) as $n\to\infty$. 
	It is well known that, for Ewens permutations $\sigma_n$ on $\SS_n$, $n\ge1$,
	\begin{equation}
		\label{eq:limcounts}
		\left(C_k(\sigma_n), k \ge 1\right) \xrightarrow d (Z_k, k \ge 1)
	\end{equation}
	in $\Z_+^\infty$, where $Z_k$ are independent and Poisson 
	distributed with means $\theta/k$, see Theorem 5.1 in \cite{ABT03}. 
	In the case of more general permutations with cycle weights, 
	Corollary 2.2 in \cite{EU14} shows that \eqref{eq:limcounts} 
	remains true with independent Poisson distributed $Z_k$ with means $\theta_k/k$ provided 
	the stability condition
	\begin{equation}
		\label{eq:stab}
		\lim_{n\to\infty} \frac{h_{n-1}}{h_n} = 1
	\end{equation}
	holds. In the same paper, it is shown that \eqref{eq:stab} 
	is satisfied for a wide range of asymptotics of $\theta_k$, 
	from sub-exponential decay to sub-exponential growth.
	
	What can be said about the limiting composition of
	short cycles themselves? Say, for fixed points (that is, $1$-cycles), 
	the invariance of \eqref{eq:dist} under relabeling suggests that, 
	conditionally on $C_1(\sigma_n)=c_1$, the set of fixed points of 
	$\sigma_n$ is distributed as a random equiprobable $c_1$-sample 
	from $[n]\vcentcolon=\{1,\ldots,n\}$ without replacement. 
	Similar reasoning can be given for cycles of any fixed length. 
	However, a rigorous description of the limiting behavior is possible 
	only within the framework of convergence of random point measures.
	Using this approach, in Section \ref{sec:main} we state and prove 
	a multivariate point processes version of \eqref{eq:limcounts}.
	Its main advantage is that it allows us to easily prove further limit theorems for different statistics of short cycles, bypassing involved combinatorial calculations and complicated asymptotic analysis. Various examples of such 
	results are given in Section \ref{sec:appl}.
	
	\section{Preliminaries and main result}
	\label{sec:main}
	
	We first introduce a metric space appropriate for describing 
	the limiting composition of cycles. For $k\ge1$, let
	\begin{equation}
		\label{eq:Xk}
		\X_k = \left\{\x=(x_1,\ldots,x_k) \in [0,1]^k \colon  \min\{x_1,\ldots,x_k\} = x_1\right\}
	\end{equation}
	and denote by $\rho_k$ the Euclidean metric on $\X_k$. 
	The last equality in \eqref{eq:Xk} is due to the fact that 
	any element of a cycle can be regarded as its 
	\lq\lq beginning\rq\rq.
	
	Consider now a multi-level space 
	$\X=\bigcup_{k=1}^\infty\X_k$ with metric given by
	\begin{equation*}
		\rho(\x_1,\x_2)=
		\begin{cases}
			\rho_k(\x_1,\x_2), & \x_1, \x_2 \in \X_k,\\
			\sqrt{\max\{k_1, k_2\}}, & \x_1 \in \X_{k_1}, \, \x_2 \in \X_{k_2}, \, k_1 \ne k_2.
		\end{cases}
	\end{equation*}
	The triangle inequality holds since
	$\sup_{\x_1,\x_2 \in \X_k} \rho_k(\x_1, \x_2) = \sqrt k$.
	The metric $\rho$ makes $\X$ a Polish space.
	Moreover, $\X$ equipped with the Borel $\sigma$-algebra $\B(\X)$ turns 
	into a measurable space, and $B\in\B(\X)$ if and only if 
	$B\cap\X_k\in\B(\X_k)$ for all $k$. This makes it possible 
	to define a measure on $\left(\X,\B(\X)\right)$ by
	\begin{equation}
		\label{eq:mu}
		\lambda(B) = \sum_{k=1}^\infty \theta_k \hs \lambda_k(B\cap\X_k), \qquad B \in \B(\X),
	\end{equation}
	where $\theta_k$ are defined in \eqref{eq:dist} and $\lambda_k$ stands for the $k$-dimensional Lebesgue measure.
	
	Let $\delta_{\x} = \mathds 1\{\x\in\cdot\}$ be the Dirac measure at $\x$. 
	We will focus on the limiting behavior of random point measures on 
	$\left(\X, \B(\X)\right)$ given by
	\begin{equation}
		\label{eq:Psin}
		\Psi_n = \sum_{k=1}^\infty\;\;
		\sideset{}{^{\ne}}\sum_{i_1,\ldots,i_k \in [n]}
		\delta_{\left(\frac{i_1}{n}, \ldots, \frac{i_k}{n} \right)}
		\mathds 1\!\left\{
		\sigma_n(i_1)=i_2, \ldots, \sigma_n(i_k)=i_1
		\right\},
	\end{equation}
	where $\sum^{\ne}$ means that the sum is taken over $k$-tuples with distinct entries. Moreover, since $\Psi_n$ is considered as a measure on $\left(\X, \B(\X)\right)$, the latter sum includes only those tuples in which the minimum element comes first.
	$\Psi_n$ carries all the 
	information about the cycle structure of $\sigma_n$.
	
	Recall that the vague topology on the space of locally finite 
	measures is generated by the integration maps 
	$\nu\mapsto\int_\X f\,\d\nu$ for all continuous functions 
	$f$ with bounded support; see, e.g., \citet[Section 3.4]{R87} 
	or \citet[Chapter 4]{K17} for a general exposition. 
	Denote by $\xrightarrow{vd}$ the distributional convergence 
	of random point measures with underlying vague topology.
	
	Let $\Psi$ denote the Poisson random measure on 
	$\left(\X, \B(\X)\right)$ with intensity measure $\lambda$ 
	given by \eqref{eq:mu}. The following theorem can be regarded as 
	a point processes extension of \eqref{eq:limcounts}.
	
	\begin{thm}
		\label{th:main}
		Let the stability condition \eqref{eq:stab} hold.
		Then $\Psi_n \xrightarrow{vd} \Psi$ as $n \to \infty$.
	\end{thm}
	
	\begin{rmk}\label{rmk:restrict}
		Let $\Psi_n^{(k)}$ and $\Psi^{(k)}$ be the restrictions to $\X_k$ of 
		$\Psi_n$ and $\Psi$, respectively.
		By the restriction property of Poisson processes (see, e.g., Theorem 5.2 in \cite{LP18}), $\Psi^{(k)}$ are independent homogeneous Poisson processes with intensities $\theta_k$. Due to the vague continuity of the restriction mapping $\mu\mapsto\mu\hspace{-3.5pt}\restriction_{\X_k}$, we also have $\Psi_n^{(k)}\xrightarrow{vd}\Psi^{(k)}$. In particular, as 
		$C_k(\sigma_n) = \Psi_n^{(k)}(\X_k)$, \eqref{eq:limcounts} directly follows from Theorem \ref{th:main}.
	\end{rmk}
	
	For the proof, we will first need an asymptotics for probabilities of cycles.
	
	\begin{lem}
		\label{lem:Pcycles}
		Let $r \ge 1$ and $\i^{(j)} = \bigl(i_1^{(j)},\ldots,i_{k_j}^{(j)}\bigr)$, 
		$j \in [r]$, be disjoint integer tuples with distinct entries such that $\i^{(j)}/n \in \X_{k_j}$. Then, under \eqref{eq:stab}, we have
		\begin{equation}
			\label{eq:Pcycles}
			\P{\sigma_n \textit{ contains the cycles } \i^{(1)},\ldots,\i^{(r)}}
			\sim \frac{\prod_{j=1}^r\theta_{k_j}}{n^{k_1+\ldots+k_r}}, \qquad n\to\infty.
		\end{equation}
	\end{lem}
	
	\begin{pf} Let $\I$ be the set of all entries of $\i^{(j)}$, $j\in[r]$, 
		and $s = \#\I = k_1 + \ldots + k_r$, 
		where $\#$ stands for the cardinality of a set. 
		By \eqref{eq:dist}, the probability in \eqref{eq:Pcycles} equals
		\begin{equation}
			\begin{aligned}
				\label{eq:Pcycles_proof}
				\sum_{\tilde\pi \in \SS_{[n]\setminus\I}} \!
				\P{
					\sigma_n = \i^{(1)} \circ \ldots \circ
					\i^{(r)} \circ \tilde\pi
				} & = \sum_{\tilde\pi \in \SS_{[n]\setminus\I}}
				\! \frac{1}{h_nn!} \prod_{k=1}^\infty \theta_k^{\#\{j\colon k_j=k\} + C_k(\tilde\pi)} \\ =
				\frac{1}{h_nn!} \prod_{j=1}^r \theta_{k_j} \cdot\!
				\sum_{\tilde\pi \in \SS_{[n]\setminus\I}}
				\prod_{k=1}^\infty \theta_k^{C_k(\tilde\pi)} & =
				\frac{h_{n-s}\hs(n-s)!}{h_nn!} \prod_{j=1}^r \theta_{k_j},
			\end{aligned}
		\end{equation}
		where the last equality follows again from \eqref{eq:dist}.
		Hence, the claim follows from \eqref{eq:stab}.\qed
	\end{pf}
	
	\begin{pot}{Theorem \ref{th:main}}
		Let $\langle$ mean either $($ or $[$, and the same applies to $\rangle$.
		By Theorem 4.18 in \cite{K17}, it suffices to prove that
		\begin{enumerate}[(i)]
			\item $\lim_{n\to\infty} \E\hs\Psi_n(B) = \E\hs\Psi(B)$ for any box
			$B = \bigtimes_{i=1}^k\langle a_j,b_j\rangle\in \X_k$, $k \ge 1$,
			\item $\lim_{n\to\infty} \P{\Psi_n(U)=0} = \P{\Psi(U)=0}$ for any finite union $U$ of boxes from possibly different levels $\X_k$.
		\end{enumerate}
		
		Let $B$ be a box with fixed $k\ge1$ and $B_{\ne}$ the set of points in $B$ with distinct coordinates. By \eqref{eq:Psin}, we have
		\begin{equation*}
			\E\hs\Psi_n(B) = \sideset{}{^{\ne}}\sum_{\left(i_1,\ldots,i_k\right)/n\in B}
			\P{\sigma_n(i_1) = i_2, \ldots, \sigma_n(i_k) = i_1}.
		\end{equation*}
		It follows from \eqref{eq:Pcycles_proof} and \eqref{eq:Pcycles} 
		with $r=1$ that all summands on the right-hand side are equal and 
		asymptotically equivalent to $\theta_k/n^k$. Thus, as $n\to\infty$,
		\[
		\E\hs\Psi_n(B) \sim \frac{\theta_k}{n^k} \cdot
		\#\left(B_{\ne}\cap(\Z^k/n)\right) \to \theta_k\hs\lambda_k(B) = \lambda(B) = \E\hs\Psi(B),
		\]
		which proves (i).
		
		We now proceed to (ii). Let us fix a finite union $U$ of boxes and denote $U_m = U \cap \X_m$. So, $U_m = \varnothing$ for $m$ greater than some $k \ge 1$ (the maximum dimension of boxes in the union),
		and $U = \bigcup_{m=1}^k U_m$.
		Let $\i_m$ and $\i_m^{(\cdot)}$ stand for integer $m$-tuples with distinct entries and $\bigcirc_{j_m=1}^{r_m} \i_m^{(j_m)}$ be the composition
		$\i_m^{(1)} \circ \ldots \circ \i_m^{(r_m)}$ of $r_m$ cycles defined by such tuples.
		For any $R\ge1$, by Bonferroni's inequality, 
		\begin{align}
			\notag
			& \P{\Psi_n(U)=0} = 1 - \PP\biggl\{\bigcup_{m=1}^k\bigcup_{\i_m/n\in U_m}
			\{\sigma_n \textit{ contains the cycle } \i_m\}\biggr\}
			\\ \label{eq:zero_prob_incl_excl}
			& \le \sum_{r_1, \ldots, r_k = 0}^{2R} (-1)^{r_1 + \ldots + r_k}
			\sideset{}{^\ast}
			\sum_{\substack{
					\i_1^{(1)}/n, \ldots, \i_1^{(r_1)}/n \in U_1, \\
					\svdots \\
					\i_k^{(1)}/n, \ldots, \i_k^{(r_k)}/n \in U_k
			}}
			\P{
				\sigma_n \textit{ contains } \bigcirc_{j_1=1}^{r_1} \i_1^{(j_1)} \circ \ldots \circ
				\bigcirc_{j_k=1}^{r_k} \i_k^{(j_k)}
			},
		\end{align}
		where $\sum^\ast$ means that the sum is taken over all unordered sets of disjoint tuples, and
		a similar lower bound holds with $2R$ replaced by $2R-1$.
		
		Due to \eqref{eq:Pcycles_proof}, all the summands in $\sum^\ast$ are equal and the right-hand side of \eqref{eq:zero_prob_incl_excl} becomes
		\begin{equation}
			\label{eq:sum}
			\sum_{r_1, \ldots, r_k = 0}^{2R} (-1)^{r_1 + \ldots + r_k}\cdot
			\frac{
				h_{n - \sum_{m=1}^k m r_m} \bigl(n - \sum_{m=1}^k m r_m\bigr)!
			}{h_n n!}
			\cdot \frac{
				\prod_{m=1}^k\theta_m^{r_m}
			}{\prod_{m=1}^k r_m!}
			\cdot S_n(r_1, \ldots, r_k),
		\end{equation}
		where
		\[S_n(r_1, \ldots, r_k) = 
		\hspace{-15pt}\sum_{\substack{
				\i_1^{(1)}/n, \ldots, \i_1^{(r_1)}/n \in \X_1, \\
				\svdots \\
				\i_k^{(1)}/n, \ldots, \i_k^{(r_k)}/n \in \X_k
		}} \hspace{-15pt}\mathds{1}\{
		\textit{all tuples are disjoint}\hs\} \prod_{m=1}^k
		\mathds{1}\bigl\{
		\i_m^{(1)}/n, \ldots, \i_m^{(r_m)}/n \in U_m
		\bigr\},\]		
		and division by $\prod_{m=1}^k r_m!$ is due to the fact that the sum in the definition of $S_n$ is taken over ordered sets of tuples.
		Note that $\frac{S_n(r_1, \ldots, r_k)}{n^{r_1 + \ldots + k r_k}}$
		can be viewed as an integral sum for
		\[
		\int_{U_1^{r_1} \times \ldots \times U_{k}^{r_k}} \mathds{1}\left\{
		\textit{all components of all } \x_{m}^{(j_m)} \textit{ are distinct}
		\right\} \;
		\prod_{m=1}^k\prod_{j_m=1}^{r_m}\d\x_m^{(j_m)} = 
		\prod_{m=1}^k \left(
		\lambda_m(U_m)
		\right)^{r_m}.
		\]
		Hence, by \eqref{eq:stab}, each summand in \eqref{eq:sum} converges as $n\to\infty$ to 
		\[(-1)^{r_1 + \ldots + r_k}\prod_{m=1}^k \frac{\left(\theta_m\lambda_m(U_m)
			\right)^{r_m}}{r_m!}.\]		
		It now follows from \eqref{eq:zero_prob_incl_excl} and a similar lower bound that
		\begin{align*}
			\sum_{r_1, \ldots, r_k = 0}^{2R-1}
			(-1)^{r_1 + \ldots + r_k}\prod_{m=1}^k \frac{\left(\theta_m\lambda_m(U_m)
				\right)^{r_m}}{r_m!}&\le\lim_{n\to\infty}\P{\Psi_n(U)=0}\\&\le\sum_{r_1, \ldots, r_k = 0}^{2R}
			(-1)^{r_1 + \ldots + r_k}\prod_{m=1}^k \frac{\left(\theta_m\lambda_m(U_m)
				\right)^{r_m}}{r_m!}.
		\end{align*}		
		Letting $R\to\infty$ finally yields
		\begin{equation*}
			\begin{aligned}
				\lim_{n\to\infty} \P{\Psi_n(U)=0} & = 
				\sum_{r_1, \ldots, r_k = 0}^{\infty} (-1)^{r_1 + \ldots + r_k}\prod_{m=1}^k \frac{\left(\theta_m\lambda_m(U_m)
					\right)^{r_m}}{r_m!} \\
				& = \exp\hs\Bigl\{ -\sum_{m=1}^k \theta_m  \lambda_m(U_m)\Bigr\} = 
				\exp\hs\{-\lambda(U)\}=\P{\Psi(U) = 0}.
			\end{aligned}
		\end{equation*}		
		This concludes the proof of (ii) and, hence, that of Theorem \ref{th:main}.\qed
	\end{pot}
	
	\section{Limit theorems for statistics of short cycles}
	\label{sec:appl}
	
	Theorem \ref{th:main} allows us to derive limiting distributions for various statistics of short cycles.
	In what follows, we will always assume that the stability condition \eqref{eq:stab} is satisfied.
	
	We first give a general result which covers the case of additive statistics.
	
	\begin{prop}
		\label{prop:add}
		Let $k\ge1$ and $f_m\colon\X_m\to[0,\infty)$, $m\in[k]$, be a family of continuous functions. Denote by $\C_m(\sigma_n)$ the set of all $m$-cycles in $\sigma_n$. Then
		\begin{equation}
			\label{eq:addstat}
			\sum_{m=1}^k\sum_{c\hs\in\hs\C_m(\sigma_n)}f_m\left(\frac cn\right)\xrightarrow d S,\qquad n\to\infty,
		\end{equation}
		where $c\in\C_m(\sigma_n)$ is understood as an integer tuple whose minimum element comes first, and
		the limiting random variable $S$ is defined by its Laplace transform
		\begin{equation}
			\label{eq:Laplace}
			\E\exp\hs\{-tS\}=
			\exp\hs\Bigl\{-\sum_{m=1}^k\theta_m\int_{\X_m}
			\left(1-\e^{-tf_m(\x)}\right)\,\d\x\Bigr\},\qquad t\ge0.
		\end{equation}
	\end{prop}
	
	\begin{pf}
		Define $f\colon\X\to[0,\infty)$ by
		$f(\x)=\sum_{m=1}^kf_m(\x)\mathds 1\{\x\in\X_m\}$.
		The left-hand side of \eqref{eq:addstat} can be written as
		$\int_\X f(\x)\,\Psi_n(\d\x)$, and $f$ is continuous with bounded support. By Theorem \ref{th:main} and Lemma 4.12 in \cite{K17}, \eqref{eq:addstat} holds with
		$S=\int_\X f(\x)\,\Psi(\d\x)$.
		The Laplace transform of $S$ is thus of the form
		\[\E\exp\hs\{-tS\}=\E\exp\hs\Bigl\{-\int_\X tf(\x)\,\Psi(\d\x)\Bigr\},\]
		which coincides with the right-hand side of \eqref{eq:Laplace} due to the form of the Laplace functional of a Poisson random measure, see, e.g., Proposition 3.6 in \cite{R87}.\qed
	\end{pf}
	
	As an example, we give a limit theorem for the sum $S_n^{(k)}$ of elements in all $k$-cycles of $\sigma_n$.
	
	\begin{prop}\label{prop:sums}~
		\begin{enumerate}[(i)]
			\item $\frac{S_n^{(k)}}n\xrightarrow d S^{(k)}$ as $n\to\infty$, where $S^{(k)}$ is defined by its Laplace transform
			\begin{equation}
				\label{eq:LaplaceSk}
				\E\exp\left\{-tS^{(k)}\right\}=\exp\hs\biggl\{\frac{\theta_k}k
				\biggl(\left(\frac{1-\e^{-t}}t\right)^k-1\biggr)\biggr\},\qquad t>0.
			\end{equation}
			\item If $k=1$, that is, for $S_n^{(1)}$ being the sum of fixed points,
			\begin{equation}
				\label{eq:cdfS1}
				\P{S^{(1)}\le x}=\e^{-\theta_1}\sum_{j=0}^{\lfloor x\rfloor}
				\frac{(-1)^j}{j!}(\theta_1(x-j))^{\frac{j}{2}}I_j\!\left(2\sqrt{\theta_1(x-j)}\right),\qquad x\ge0,
			\end{equation}
			where $I_j$ is the modified Bessel function of the first kind, see, e.g., \S\hs10.25\hs(ii) in \cite{NIST}.
		\end{enumerate}	
	\end{prop}
	
	\begin{pf}
		For $\x\in\X_k$, let $f_k(\x)$ denote the sum of all its components. By Proposition \ref{prop:add}, (i) holds with
		\[\E\exp\left\{-tS^{(k)}\right\}=\exp\hs\Bigl\{-\theta_k\int_{\X_k}
		\left(1-\e^{-tf_k(\x)}\right)\,\d\x\Bigr\},\]
		where the integral on the right-hand side, due to symmetricity of $f_k$, equals
		\[\frac 1k\int_{[0,1]^k}
		\left(1-\e^{-tf_k(\x)}\right)\,\d\x=
		\frac 1k\Bigl(1-\Bigl(\int_0^1\e^{-tx}\,\d x\Bigr)^k\,\Bigr).\]
		This yields \eqref{eq:LaplaceSk}.		
				
		To prove \eqref{eq:cdfS1}, we first note that
		\[\E\exp\left\{-tS^{(1)}\right\}=
		\E\int_0^\infty t\e^{-tx}\,\mathds 1\{S^{(1)}\le x\}\,\d x
		=t\int_0^\infty \e^{-tx}\,\PP\{S^{(1)}\le x\}\,\d x,\]
		cf.~\cite{L20}. Hence, $\P{S^{(1)}\le x}$ is the inverse Laplace transform of the right-hand side in \eqref{eq:LaplaceSk} for $k=1$ multiplied by $t^{-1}$, that is, in expanded form, of the function
		\begin{equation}
		\label{eq:Ginf}
		G(t)=\sum_{j=0}^\infty\frac{(-\theta_1)^j}{j!}\e^{-\theta_1}\hs t^{-j-1}\e^{\frac{\theta_1}t}\e^{-jt},\qquad t>0.
		\end{equation}
		By the time shifting property and \citet[eq.~(5.5.35)]{E54},
		the inverse Laplace transform of the $j$-th summand $G_j(t)$ in \eqref{eq:Ginf} equals
		\[F_j(x)=\frac{(-1)^j}{j!}e^{-\theta_1}\hs(\theta_1(x-j))^{\frac j2}I_j\!\left(2\sqrt{\theta_1(x-j)}\right)\mathds 1\{x\ge j\}.\]
		This means that, for all $t>0$ and $R\in\mathbb N$,
		\begin{equation}
		\label{eq:R}
		\int_0^\infty\!\e^{-tx}\sum_{j=0}^R F_j(x)\,\d x=\sum_{j=0}^R G_j(t).
		\end{equation}
		Since $I_j(x)$ decreases in $j$ and increases in $x$ (see, e.g., \S\hs10.37 in \cite{NIST}), and $I_0(x)\sim\frac{\e^x}{\sqrt{2\pi x}}$ as $x\to\infty$ (ibid., \S\hs10.30\hs(ii)), we have
		\[\Bigl|\sum_{j=0}^R F_j(x)\Bigr|\le
		e^{-\theta_1}I_0\!\left(2\sqrt{\theta_1x}\right)
		\sum_{j=0}^\infty\frac 1{j!}(\theta_1x)^{\frac j2}\sim
		e^{-\theta_1}\frac{\e^{3\sqrt{\theta_1x}}}{2\pi^{\frac 12}(\theta_1x)^{\frac14}},\qquad x\to\infty,\]
		which makes it possible to apply dominated convergence to \eqref{eq:R}.
		Hence, the inverse Laplace transform of $G$ is $\sum_{j=0}^\infty F_j$,
		which yields \eqref{eq:cdfS1}.\qed
	\end{pf}
	
	We now turn to some examples of non-additive statistics. 
	For $k\ge2$, let us call the range of a cycle the difference between its 
	maximum and minimum elements and denote by $r_n^{(k)}$ (resp., $R_n^{(k)}$) the minimum (maximum) range among all $k$-cycles in $\sigma_n$.
	If there are no $k$-cycles, we set $r_n^{(k)}=n$ and $R_n^{(k)}=0$.
	
	\begin{prop}
		\label{prop:ranges}
		$\frac{r_n^{(k)}}n\xrightarrow d r^{(k)}$ 
		and $\frac{R_n^{(k)}}n\xrightarrow d R^{(k)}$ as $n\to\infty$, where
		$r^{(k)}$ {\rm(}resp., $R^{(k)}${\rm)}, $k\ge2$, have CDF's of the form
		\begin{gather}
			\label{eq:range1}
			\P{r^{(k)}\le x}=1-
			\exp\hs\Bigl\{-\frac{\theta_k}k\left(kx^{k-1}-(k-1)x^k\right)\Bigr\},\\
			\P{R^{(k)}\le x}=
			\exp\hs\Bigl\{\frac{\theta_k}k\left(kx^{k-1}-(k-1)x^k-1\right)\Bigr\},
			\label{eq:range2}
		\end{gather}
		as $x\in[0,1)$ and $0$ {\rm(}resp., $1${\rm)} to the left
		{\rm(}right\hs\hs{\rm)}.
	\end{prop}
	
	\begin{pf} 	For $\x\in\X_k$, let $f_k(\x)$ denote the difference between its maximum and minimum components. It follows from the interpretation of vague convergence in Proposition 3.13 of \cite{R87} that the function which maps a finite point measure $\mu$ on $\X_k$ into $\min_{\mu\{\x\}\ge1}f_k(\x)$ is vaguely continuous. Hence, by Remark \ref{rmk:restrict} and continuous mapping theorem, $\frac{r_n^{(k)}}n\xrightarrow d r^{(k)}$ holds with $r^{(k)}=\min_{\Psi^{(k)}\{\x\}\ge1}f_k(\x)$, and the CDF of $r^{(k)}$ is of the form
		\begin{equation}
			\begin{aligned}
				\label{eq:range_prob}
				\P{r^{(k)}\le x}&=1-\P{\Psi^{(k)}\{\x\in\X_k\colon f_k(\x)\le x\}=0}
				\\&=1-\exp\hs\bigl\{-\theta_k\hs\lambda_k
				\{\x\in\X_k\colon f_k(\x)\le x\}\bigr\}.
			\end{aligned}
		\end{equation}
		The value of the Lebesgue measure on the right-hand side is
		\[\int_0^1\biggl(\int_{x_1}^{\min\{x_1+x,1\}}\d x_2\ldots\d x_k\biggr)\d x_1=x^{k-1}(1-x)+\frac{x^k}k,\]
		which together with \eqref{eq:range_prob} yields \eqref{eq:range1}.
		The proof of \eqref{eq:range2} is similar.\qed
	\end{pf}
	
	As a final example, we consider some statistics of fixed points. Let $m_n$ be the minimum fixed point ($n+1$ in case there is none), $M_n$ the maximum one ($0$ in that case), $\delta_n$ the minimum spacing between fixed points (the two extreme spacings of lengths $m_n$ and $n+1-M_n$ are also taken into account), and $\Delta_n$ the maximum one.
	
	\begin{prop}\label{prop:fixed_points_stats}
		$\left(\frac{m_n}n, \frac{M_n}n, \frac{\delta_n}n, \frac{\Delta_n}n\right)
		\xrightarrow d(m,M,\delta,\Delta)$ as $n\to\infty$ with
		\begin{gather*}
			\P{m\le x}=1-\exp\hs\{-\theta_1x\},\qquad
			\P{M\le x}=\exp\hs\{\theta_1(x-1)\},\qquad x\in[0,1),\\
			\delta\overset d=\frac{X_{\nu+1}}{(\nu+1)\sum_{i=1}^{\nu+1}X_i},\qquad
			\Delta\overset d=\frac{\sum_{i=1}^{\nu+1}\frac{X_i}i}{\sum_{i=1}^{\nu+1}X_i},
		\end{gather*}
		where $\nu$ is Poisson distributed with parameter $\theta_1$, $X_i$ are exponentially distributed with unit mean, and all these are independent.
	\end{prop}
	
	\begin{pf}
			As in the proof of Proposition \ref{prop:ranges}, the convergence takes place with $m$ being the leftmost point of $\Psi^{(1)}$, $M$ the rightmost one, $\delta$ the minimum spacing, and $\Delta$ the maximum one. All that remains is to derive the corresponding CDF's and distributional equalities.			
			For $m$ and $M$, this follows from
			\begin{gather*}
			\P{m\le x} = 1-\P{\Psi^{(1)}[0, x] = 0} = 1-\exp\hs\{-\theta_1x\},\\
			\P{M\le x} = \P{\Psi^{(1)}(x, 1] = 0} = \exp\hs\{-\theta_1(1-x)\}.
			\end{gather*}			
					
			We now turn to the equalities for $\delta$ and $\Delta$.
			Since $\Psi^{(1)}$ is a homogeneous Poisson process with intensity $\theta_1$, the conditional distribution of $\delta$ (resp., $\Delta$) given $\Psi^{(1)}(\X_1) = r$ coincides with that
			of the minimum (maximum) spacing $d_r$ ($D_r$) for $r$ independent random variables, uniformly distributed on $[0,1]$, see, e.g., Proposition 3.8 in \cite{LP18}.
			
			It follows from the theory of uniform spacings (see, e.g., \cite{H80}, p.~625) that
			\begin{equation}
				\label{eq:Uspac}
				d_r\overset d=\frac{\min Y_i}
				{\sum_{i=1}^{r+1}Y_i},
				\qquad D_r\overset d=\frac{\max Y_i}
				{\sum_{i=1}^{r+1}Y_i},
			\end{equation}
			where $Y_1,\ldots,Y_{r+1}$ are independent $\mathsf{Exp}(1)$-distributed random variables. Let $Y_{(i)}$, $i\in[r+1]$, be their order statistics, $Y_{(0)}=0$, and $\tau_i=Y_{(i)}-Y_{(i-1)}$. By Sukhatme--R\'enyi decomposition (see, e.g., Theorem 4.6.1 in \cite{ABN08}), $\tau_i$ are independent and $\mathsf{Exp}(r-i+2)$-distributed
			with
			\[\min Y_i=\tau_1,\qquad
			\max Y_i=\sum_{i=1}^{r+1}\tau_i,\qquad
			\sum_{i=1}^{r+1}Y_i=\sum_{i=1}^{r+1}(r-i+2)\tau_i.\]
			Denoting $X_{r-i+2}=(r-i+2)\tau_i \sim \mathsf{Exp}(1)$, we can rewrite \eqref{eq:Uspac} as
			\[
			d_r \overset d= \frac{X_{r+1}}{(r+1) \sum_{i=1}^{r+1} X_i}, \qquad
			D_r \overset d= \frac{\sum_{i=1}^{r+1}\frac{X_i}i}
			{\sum_{i=1}^{r+1} X_i}.
			\]
			Since this holds for any $r\ge0$, the claim follows by deconditioning.
	\end{pf}
		
	\bibliographystyle{elsarticle-harv}
	\bibliography{Cycles}
	
\end{document}